\documentclass[reqno, 11pt] {amsproc}

\usepackage[cp1251] {inputenc}
\usepackage[english] {babel}
\usepackage{amsmath, amssymb, babel, amsthm}
\newtheorem{theorem} {Theorem}
\newtheorem{lemma} {Lemma}
\newtheorem{proposition} {Statement}
\newtheorem{remark} {Note}

\title[Forward and inverse problems for a finite Krein-Stieltjes string.]
{Forward and inverse problems for a finite Krein-Stieltjes string.
Approximation of constant density by point  masses.}

\author{A.S. Mikhaylov}
\address{St. Petersburg branch of the Mathematical Institute. V.A. Steklova RAS, 191023,
    nab. Fontanka 27, St. Petersburg, Russia and St. Petersburg
    State University, 199034, Universitetskaya emb. 7/9,
    St. Petersburg, Russia.} \email{mikhaylov@pdmi.ras.ru}

\author{V.S. Mikhaylov}
\address{St. Petersburg branch of the Mathematical Institute. V.A. Steklova RAS, 191023,
    nab. Fontanka 27, St. Petersburg, Russia}
\email{ftvsm78@gmail.com}


\date{March 12, 2020}
\begin{document}
\maketitle
\begin{abstract}
We consider a dynamic inverse problem for a dynamical system which
describes the  propagation of waves in a Krein string. The problem
is reduced to an integral equation and an important special case
is considered when the string density is determined by a finite
number of point masses distributed over the interval. We derive an
equation of Krein type, with the help of which the string density
is restored. We also consider the approximation of constant
density by point masses uniformly distributed over the interval
and the effect of the appearance of a finite wave propagation
velocity in the dynamical system.
\end{abstract}

\section{Introduction}
Let a non-decreasing bounded function be given on the interval $
(0, l) $ $ M (x) $: $ 0 \leqslant M (0) <M (x) <M (l) $. By the
value of the function $ M $ at the point $ x $ we mean the mass of
the string segment $ (0, x) $. If $ M $ is a differentiable
function, then $ \rho (x) = M '(x) $ is the string density.
Following \cite{DMcK, KK}, we introduce a domain consisting of
\emph{continued functions}
\begin{multline*}
    D_M:=\Bigl\{\left[u(x),u'_-(0),u'_+(l)\right]\,\bigl|\quad
    u(x)=a+bx-\int_0^x(x-s)g(s)\,dM(s);\\
    u'_-(0)=b;\quad u'_+(l)=b-\int_0^lg(s)\,dM(s) \Bigr\},
\end{multline*}
 where $ u '_- (0), u' _ + (l) $ are left and right derivatives, $ g $ is
 $ M $-summable function. Define \emph{generalized differential Krein operation}
 $ l_M [u] $ on $ D_M $ by the rule
\begin{equation*}
    \label{L_oper} l_M[u]=g(x),\quad \text{$M-$almost everywhere}.
\end{equation*}
   Note that if the function $ M $ is differentiable and $ M '(x)> 0 $, then
    $l_m[y] = -\frac{1}{M'(x)} \frac{d^2 y(x)}{dx^2}$.

We fix $ T> 0 $ and consider the following initial-boundary value
problem:
\begin{equation}
    \label{L_eq}
    \begin{cases}
        u_{tt}(x,t)+l_M[u]=0,\quad 0<x<l,\, 0<t<T,\\
        u(x,0)=0,\, u_t(x,0)=0,\quad 0\leqslant x\leqslant l,\\
        u(0,t)=f(t),\,  u(x,l)=0\quad 0\leqslant t\leqslant T,
    \end{cases}
\end{equation}
where $ f \in L_2 (0, T) $ \, - \, is a \emph{boundary control}.

Note that the authors did not find references in the literature to
dynamic inverse problems for dynamical systems (\ref{L_eq}),
except for the cases when $ M '(x) = \rho(x) \geqslant \delta> 0
$, $ 0 \leqslant x \leqslant l $, then for $ \rho \in C ^ 1 $ the
system (\ref{L_eq}) corresponds to the initial-boundary value
problem for the wave equation \cite{BM, BL}. The reason that such
problems have not yet been investigated is the fact that the wave
propagation velocity in systems with a `` bad '' density is not
finite and usual methods, such as the Boundary Control method
\cite{AM, B07, B17, BM} are not directly applicable, and should be
modified.

In this paper, we want to present a general strategy for studying
the dynamic inverse problem for the system (\ref{L_eq}) using the
Boundary Control method \cite{B07, B08, B17} based on a special
case - the finite Krein-Stieltjes string.

In the second section, for an arbitrary mass distribution, we
derive the integral equation equivalent to (\ref{L_eq}). In the
third section, following \cite{MM4}, we consider direct and
inverse problems in an important special case: the finite
Krein–Stieltjes string. In the last section, we consider an
example where the unit density is approximated by point masses
uniformly distributed over the interval.

\section{The direct problem for Krein string}
It is useful to rewrite (\ref{L_eq}) in integral form. Twice
integrating the equation in (\ref{L_eq}), we get that
\begin{equation*}
u(x,t)=-\int_0^t(t-s)l_M[u(\cdot,s)]\,ds.
\end{equation*}
Multiplication by $ x $ and integration gives
\begin{equation}
\int_0^yu(x,t)x\,dM(x)\label{Eq1}
=-\int_0^y\int_0^t(t-s)xl_M[u(x,s)]\,ds\,dM(x).
\end{equation}
Using Lagrange identities \cite{DMcK}
\begin{equation*}
\int\limits_0^l\left(l[u](x)v(x) - u(x)l[v](x) \right)\,dM(x) =u(x)v'(x)-u'(x)v(x)\Big|_{x=0-}^{x=l+}
\end{equation*}
 with the function $ v (x) = x $ we arrive at
\begin{equation*}
\int_0^l xl_M[u(x,s)]\,dM(x)=u(x,t)-u'(x,t)\Big|_{x=0-}^{x=l_+},
\end{equation*}
and thus (\ref{Eq1}) reduces to
\begin{multline*}
\int_0^yu(x,t)x\,dM(x)= -\int\limits_0^t (t-s)(u(x,s)-u'(x,s)x)\Big|_{x=0}^{x=y+} ds         \\
=-\int_{0}^t(t-s)\left(u(y,s)-u(0,s)-yu'(y+,s)\right)\,ds.
\end{multline*}
Integrating     the last identity with respect to the variable $ y $ from $ 0 $ to $ z $ we get
\begin{multline*}
\int_0^z\int_0^yu(x,t)x\,dM(x)\,dy=-\int_{0}^t(t-s)\int_0^z\left(u(y,s)-u(0,s)\right)\,dy\,ds+\\
+\int_{0}^t(t-s)\int_0^z\left(yu'(y+,s)\right)\,ds= -2\int_0^t(t-s)\int_0^z u(y,s)\,dy\,ds +\\
+ z\int_0^t(t-s)\left[u(0,s)+u(s,z)\right]\,ds.
\end{multline*}
All calculations lead to the statement:
\begin{lemma}
The initial boundary value problem (\ref{L_eq}) is equivalent to the following integral equation:
\begin{multline}
\label{Krein_int}\int_0^x
z(x-z)u(z,t)\,dM(z)=x\int_0^t(t-s)\left[f(s)+u(x,s)\right]\,ds\\
-2\int_0^t(t-s)\int_0^x u(y,s)\,dy\,ds,\quad u(x,l)=0.
\end{multline}
\end{lemma}
Solvability of (\ref{Krein_int}) for an arbitrary distribution
mass $ M $ and analysis of the \emph{control operator} defined by the rule
\begin{equation*}
W^T: f(t)=u(0,t)\mapsto u^f(\cdot,T),
\end{equation*}
are an important step in the procedure for solving dynamic inverse
problems, see \cite{B07, B08, B17, BM, MM, MM3}. Now we consider a
special case of mass distribution, namely the Krein-Stieltjes
string.

\section{Special case: finite Krein-Stieltjes string. Forward and inverse problems.}

When studying the spectral problems \cite{A, DMcK, KK}, the "bad"
mass  is approximated by masses that have the form of steps, which
means that the corresponding densities have a simple structure of
point masses distributed over the interval. That is why the
results of this section, more fully described in \cite{MM4}, can
be used as the first step to solving the dynamic inverse problem
for a general string.

We assume that the string is a Krein-Stieltjes string, i.e. the
mass $ M $ is a piecewise constant function. Let
$0=x_0<x_1<x_2<\ldots<x_{N-1}<x_N=l$, $m_i>0,$ $i=1,\ldots,N-1$,
$l_i=x_i-x_{i-1}$, $i=1,\ldots N,$ and the density $ dM $ have the
form $dM(x)=\sum_{i=1}^{N-1}m_i\delta(x-x_i)$. In this case, we
find that
\begin{align*}
&l_M[u](x)=\frac{u'(x+0)-u'(x-0)}{m_x},\\
&m_x:=\begin{cases} M(+0)-M(0),\quad x=0,\\
M(x+0)-M(x-0),\quad 0<x<l,\\
M(l)-M(l-0),\quad x=l,
\end{cases}
\end{align*}
and, therefore, $ u '(x) = 0 $, $ x \in (x_ {i-1}, x_i) $ and
introducing the notation $ u_i (t) = u (x_i, t) $, we see that the
initial-boundary value problem (\ref{L_eq}) is equivalent to the
initial-boundary value problem for a vector function (we keep the
same notation) $u(t)=\left(u_1(t),\ldots,u_{N-1}(t)\right)$:
\begin{equation}
\label{DynSyst}
\begin{cases}
Mu_{tt}=-Au+\widetilde{f},\quad t> 0,\\
u(0)=0,\, u_t(0)=0,
\end{cases}
\end{equation}
where
\begin{eqnarray*}
A=\begin{pmatrix} b_1 & a_1&0&   \ldots & 0\\
a_1& b_2 & a_2& \ldots & 0\\
\ldots & \ldots& \ldots & \ldots \\
0 & 0& \ldots&  b_{N-2} &a_{N-2}\\
0 & 0& \ldots& a_{N-2} &b_{N-1}
\end{pmatrix},\\
M=\begin{pmatrix} m_1 & 0&  \ldots & 0\\
0& m_2 &  \ldots & 0\\
\ldots & \ldots& \ldots & \ldots \\
0 & \ldots& 0 & m_{N-1}
\end{pmatrix}, \,
\widetilde f=\begin{pmatrix} \frac{f}{l_1}\\0\\
\ldots\\0\end{pmatrix},
\end{eqnarray*}
and the elements of the matrix $ A $ are defined as
\begin{equation}
\label{coeff} a_i=\frac{1}{l_{i+1}},\quad
b_i=-\frac{l_i+l_{i+1}}{l_il_{i+1}}.
\end{equation}
We note that in this case $A<0$.

The solution of the initial-boundary value problem (\ref{L_eq})
is defined as a continuous function according to the following
rule:
\begin{equation}\label{con_v_f}
u^f(x,t)=\begin{cases}  f(t),\quad x=0,\\
u_i(t), \quad x=x_i,\quad i=1,\ldots,N-1,\\
0,\quad x=l,\\
\text{affine function},\quad x\in (x_{i-1},x_i).
\end{cases}
\end{equation}
The solution to (\ref{DynSyst}) is denoted by $ u^ f (t) $. In
what follows we  fix  $T>0$.

\emph{Dynamic inverse problem} consists in restoring the string,
that is, the coefficients $ l_N $, $ m_i $, $ l_i $, $ i = 1,
\ldots, N-1 $ or, equivalently, the elements of the matrix $ A $
and $ M $, from the \emph{response operator} $ R^ T $ defined by
the rule
\begin{equation}
\label{Resp} \left(R^Tf\right)(t):=u^f_1(t),\quad 0<t<T.
\end{equation}

Consider the following Cauchy problem for a difference equation
with a spectral parameter $\lambda$:
\begin{equation*}
\begin{cases}
a_n\phi_{n+1}+a_{n-1}\phi_{n-1}+b_n\phi_n=\lambda m_n\phi_n,\\
\phi_0=0,\phi_1=1,
\end{cases}
\end{equation*}
$n = 1, \ldots, N-1$. Sequentially changing the value of the index
$ n $ from $ 1 $ to $ N-1 $,  as the solutions we obtain the set
of polynomials (of the parameter $ \lambda $) $ \{0,1,
\phi_2(\lambda), \ldots, \phi_N (\lambda) \} $. Let $ \lambda_1,
\ldots, \lambda_{N-1} $ denote the roots of the equation $ \phi_N
(\lambda) = 0 $. After introduction of vectors and coefficients
\begin{equation*}
\varphi(\lambda)=\begin{pmatrix} \phi_1(\lambda)\\
\cdot\\ \cdot\\ \phi_{N-1}(\lambda)\end{pmatrix},\quad\varphi^k=\begin{pmatrix} \phi_1(\lambda_k)\\
\cdot\\ \cdot\\ \phi_{N-1}(\lambda_k)\end{pmatrix},\quad
\omega_k=\left(M\varphi^k,\varphi^k\right),
\end{equation*}
we denote by \emph{spectral data} and \emph{spectral function}
the following objects:
\begin{equation}
\label{measure}
\left\{\lambda_i,\omega_i\right\}_{i=1}^{N-1},\quad
\mu(\lambda)=\sum_{\{k:\lambda_k<\lambda\}}\frac{1}{\omega_k}.
\end{equation}
\emph{The external space} of the system (\ref{L_eq}),
(\ref{DynSyst}),  or the space of controls, we denote
$\mathcal{F}^T:=L_2(0,T)$.
\begin{lemma}\label{Spec_lem}
The solution to the system (\ref{DynSyst}) allows spectral
representation
\begin{eqnarray}
\label{Sol_spec_repr}
u^f(t)=\frac{1}{l_1}\int_{-\infty}^\infty\int_0^t\frac{\sin{\sqrt{|\lambda_k|}(t-\tau)}}{\sqrt{|\lambda_k|}}f(\tau)\,d\tau\varphi(\lambda)\,d\mu(\lambda).\notag
\end{eqnarray}
The response operator $ R^T: \mathcal{F}^T \mapsto \mathcal{F}^T $
has the form
\begin{equation*}
\left(R^Tf\right)(t)=\int_0^t r(t-s)f(s)\,ds,
\end{equation*}
where
\begin{equation}
\label{resp_func}
r(t)=\frac{1}{l_1}\sum_{k=1}^{N-1}\frac{\sin{\sqrt{|\lambda_k|}t}}{\sqrt{|\lambda_k|}\omega_k}
\end{equation}
is called a \emph{response function}.
\end{lemma}
\emph {Inner space}, , i.e. the space of  states of system
(\ref{DynSyst}) is $ \mathcal{H}^N: = \mathbb{R}^{N-1}$. For $ T >
0 $, $ u^f(T) \in \mathcal{H}^N $. The metric in  $ \mathcal{H}^N
$ is determined by the formula
\begin{equation*}
(a,b)_{\mathcal{H}^T}=\left(Ma,b\right)_{\mathbb{R}^{N-1}},\quad
a,b\in\mathcal{H}^N.
\end{equation*}
The boundary control problems for some wave equations with
nonsmooth densities were studied in \cite{AA}. \emph{Control
operator} $W^T: \mathcal{F}^T \mapsto\mathcal{H}^N $ is determined
by the rule:
\begin{equation}
\label{Con2} W^Tf=u^f(T).
\end{equation}
We  introduce the space
\begin{equation*}
\mathcal{F}^T_1=\overline{\operatorname{span}\left\{\frac{\sin{\sqrt{|\lambda_k|}(T-t)}}{\sqrt{|\lambda_k|}}\right\}}.
\end{equation*}
The following lemma answers the question of boundary
controllability of (\ref{DynSyst}):
\begin{lemma}
\label{LemmaCont} The operator $ {W}^T $ is an isomorphism between
$\mathcal{F}^T_1$ and $\mathcal{H}^N$.
\end{lemma}
\emph{Connecting operator} $ C^T: \mathcal{F}^T \mapsto
\mathcal{F}^T $ is defined as $ C^T: = \left(W^T \right)^*W^T $,
therefore, by the definition, for $ f, g \in \mathcal{F}^T $ we
have that
\begin{equation*}
\left(C^Tf,g\right)_{\mathcal{F}^T}=\left(u^f(T),u^g(T)\right)_{\mathcal{H}^N}=\left(W^Tf,W^Tg\right)_{\mathcal{H}^N}.
\end{equation*}
In the Boundary control method, it is important that the operator
$ C ^ T $ can be expressed in terms of inverse data:
\begin{theorem}
The connecting operator can be represented in terms of dynamic
inverse data
\begin{equation*}
\left(C^Tf\right)(t)=\frac{1}{l_1}\int_0^T\int_{|t-s|}^{2T-s-t}r(\tau)\,d\tau
f(s)\,ds.
\end{equation*}
\end{theorem}
\begin{remark}
The relation $ \mathcal{F}^T_1 = C^T \mathcal{F}^T $ is valid.
It means that $\mathcal{F}^T_1 $  is determined by the inverse
data.
\end{remark}

In \cite{MM4}, the authors proposed three methods for string
recovery  from knowledge of the response operator $ R^T $ for an
arbitrary $T> 0 $. It was shown that it is enough to know the
function $ r $ and its derivatives only at zero. Here we give only
one reconstruction method based on the Krein equation.

By $ f_k^T \in \mathcal{F}^T_1 $ we denote the controls that bring
the system (\ref{DynSyst}) to the prescribed \emph{special
states}:
\begin{equation*}
d_k\in \mathcal{H}^N,
\,d_k=\left(0,.,\frac{1}{m_k},\ldots,0\right),\, k=1,\ldots,N-1.
\end{equation*}
In other words,  $ f^T_k \in \mathcal{F}^T_1$ are solutions of the
equations $ W^Tf_k^T = d_k $, $ k = 1 \, \ldots, N-1 $. According
to the lemma \ref{LemmaCont} we know that such controls exist and
are unique in  $\mathcal{F}^T_1$.  It is important for us that
they can be found as solutions of the Krein equation.
\begin{theorem}
\label{Theor_Krein_eq} The control $ f_1^ T $ can be found as
solution of the following equation:
\begin{equation*}
\left(C^Tf_1^T\right)(t)=r(T-t),\quad 0<t<T.
\end{equation*}
The controls $ f^T_k $, $ k = 1, \ldots, N-1 $ satisfy the system:
\begin{equation*}
\begin{cases}
m_1\left(C^Tf_1^T\right)''=b_{1}C^Tf^T_{1}+a_1C^Tf^T_2,\\
\cdot\cdot\cdot\cdot\cdot\cdot\cdot\cdot\cdot \\
m_k\left(C^Tf_k^T\right)''=a_{k-1}C^Tf^T_{k-1}+b_kC^Tf^T_k+a_kC^Tf^T_{k+1},\\
\cdot\cdot\cdot\cdot\cdot\cdot\cdot\cdot\cdot\\
m_{N-1}\left(C^Tf_{N-1}^T\right)''=a_{N-2}C^Tf^T_{N-2}+b_{N-1}f^T_{N-1},
\end{cases}
\end{equation*}
where $ a_0 = a_N = 0. $ String parameters can be restored as
\begin{eqnarray*}
m_k=\frac{1}{\left(C^Tf_k^T,f_k^T\right)_{\mathcal{F}^T}},\,\,
b_k=m_k^2\left(\left(C^Tf_k^T\right)'',f_k^T\right)_{\mathcal{F}^T},\\
a_k=-b_k-a_{k-1},\,\,
l_1=-\frac{\|f_1^T\|^2_{L_2(0,T)}}{\left(f_1^T\right)'(T)},\
\end{eqnarray*}
where $k=1,\ldots,N-1$.
\end{theorem}

\begin{remark}
In one-dimensional discrete and continuous dynamical systems,
\cite{B07, B08, BM, MM2, MM3, MM4} (analogue) of the control
operator $ W^T $ (\ref{Con2}) is an isomorphism between the
control space and the state space, which corresponds to exact
boundary controllability of the corresponding system. However, in
deriving the Theorems \ref{Theor_Krein_eq} we used a weaker type
of controllability, namely \emph{spectral controllability.}
\end{remark}

In our opinion, the case of the Krein-Stieltjes string with an
infinite number of point masses should be considered in the
framework of the approach described above, and not in the
framework of the general approach, when the density is
approximated by a finite number of point masses distributed over
the interval.

\section{The finite and infinite speed of wave propagation. \\ Constant density.}

In this section we consider the case when a string on the
interval $ (0,1) $ with a density $ \rho(x) = M'(x) = 1 $ is
approximated by uniformly distributed point masses. For this case,
the initial-boundary value problem (\ref{L_eq}) has the form:
\begin{equation}
\label{eq_unit}
\begin{cases}
v_{tt}(x,t)-v_{xx}=0,\quad 0<x<l,\, 0<t<T,\\
v(x,0)=0,\, v_t(x,0)=0,\quad 0\leqslant x\leqslant l,\\
v(0,t)=f(t),\,  u(x,l)=0\quad 0\leqslant t\leqslant T,
\end{cases}
\end{equation}

We emphasize that in this system the wave propagation velocity is
finite and equal to one. As an approximation, we use the identical
point masses and lengths: $ l_i = 1 / N $, $ m_i = 1 / N $, $ i =
1, \ldots, N-1 $ and, as follows from the lemma \ref{Spec_lem} ,
the wave propagation velocity in an approximate system is
infinite. Our goal is to obtain the effect of the appearance of a
finite wave propagation velocity at $ N \to \infty. $

The matrices $ A $, $ M $ in (\ref{DynSyst}) have forms:
\begin{equation}
A=N\begin{pmatrix} -2 & 1&   \ldots & 0\\
1& -2 &  \ldots & 0\\
\ldots & \ldots& \ldots &1 \\
0 &  \ldots & 1 &-2
\end{pmatrix},\quad 
M=\frac{1}{N}\begin{pmatrix} 1 & 0&  \ldots & 0\\
0& 1 &  \ldots & 0\\
\ldots & \ldots& \ldots & \ldots \\
0 & \ldots& 0 & 1
\end{pmatrix}.\label{MN}
\end{equation}

The characteristic function for the matrix $ A_M= M^{-1/2} A M^{-1/2} $ is calculated
by the formula $ \chi_{N-1} (x) = \det (A_M- \lambda I) $, where
$\lambda $ and $ x $ are related by the equality $ \lambda = -2N^2
(1 + x) $ and the index $ N $ corresponds to the size of the
matrix $ A_M \in \mathbb{M}^{N-1} $. From the recurrence relation
$ \chi_m (x) = 2x\chi_{m-1} (x) - \chi_{m-2} (x) $ for $ m = N-1,
N-2, \dots $ and that $ \chi_1 = 2x $, $ \chi_2 = 4x ^ 2-1 $ we
get that $ \chi_m (x) = U_m (x) $ where $ U_m $ \, - \ is a
Chebyshev polynomial of the second kind. Therefore, the
eigenvalues of $ A_M $ are $ \lambda_k = -2N^2 (1 + x_k) $, where
$ x_k = \cos\left(\frac{k \pi}{N} \right) $, $ k = 1, \ldots, N-1
$. Putting the first component of the eigenvector equal to unity,
we obtain
\begin{equation*}
\lambda_k=-4N^2\cos^2\left(\frac{k\pi}{2N}\right),\quad
b_k=\left(1,U_1(-x_k),\ldots,U_{N-2}(-x_k)\right),
\end{equation*}
for $ k = 1, \ldots, N-1 $ are eigenvalues and (non-normalized)
eigenvectors of $ A_M $. We calculate the norm of these
eigenvectors using the orthogonality property
\begin{equation}\label{ort_ch}
\sum\limits_{k=1}\limits^{N-1} U_i(x_k) U_j(x_k)(1-x_k^2)=\begin{cases}0,\quad i\not=j,\\
\frac{N}{2}, \quad i=j
\end{cases}
\end{equation}
from which we infer that
$|b_k|^2=\frac{N}{2\sin^2(\frac{k\pi}{N})}$.

Then from Lemma \ref{Spec_lem} for $ f (t) = \delta (t) $ we get
spectral representation for $ j = 1, \ldots, N-1 $:
$$
u^\delta_j(t)=(-1)^{j-1}\sum\limits_{k=1}\limits^{2N-1}\sin(2Nx_kt)(1-x_k^2)U_{2j-1}(x_k),
$$
where $x_k=\cos\left(\frac{k\pi}{2N}\right)$. The last expression
can be converted to integral form. To do this, one need to expand
$ \sin(2Nxt) $ in a series of $ U_j (x) $ on the segment $ [- 1,1]
$, and then using (\ref{ort_ch}) one can get the following
expression for $ u ^\delta $:
\begin{equation*}
u^\delta_j(t)=(-1)^{j-1} N\frac{2}{\pi} \int\limits_{-1}^{1}
\sin(2Nxt) U_{2j-1}(x)\sqrt{1-x^2}dx,\quad  j = 1, \ldots, N-1.
\end{equation*}
 The last integral can be calculated using the Bessel function: the integrand is written as
$$
\frac{1}{2} \frac{\sin(2Nxt)}{\sqrt{1-x^2}}
(T_{2j-1}(x)-T_{2j+1}(x)),
$$
and using the formula
$$
\int\limits_0^1 T_{2n+1}(x)\sin( ax )\frac{dx}{\sqrt{1-x^2}} =
(-1)^n\frac{\pi}{2} J_{2n+1}(a)
$$
from \cite{GR} (formula (7.355) on page 850) we come to
\begin{proposition}
The solution of the initial-boundary value problem (\ref{L_eq})
for the  Krein-Stieltjes string with uniformly distributed
identical point masses has the form (\ref{con_v_f}), where the
components of the vector $ u^\delta (t) $ are defined as
\begin{equation}\label{u_1^delta}
u^\delta_j(t)=\frac{2j}{t}J_{2j}(2Nt) = N
(J_{2j-1}(2Nt)+J_{2j+1}(2Nt))
\end{equation}
for $j=1,\ldots,N-1$,  where $J_n$\,--\, Bessel functions.
\end{proposition}
We are interested to see what happens in the limit when $N \to \infty $:
\begin{proposition}\label{p2}
The response function $ r_N (t) $ for a dynamic system
(\ref{DynSyst})  with matrices defined as (\ref{MN}) (for
the Krein-Stieltjes strings with uniformly distributed identical point masses)
converges in the sense of distributions to 
$$
r_N(t)\to  \delta(t) 
\quad\text{ for }N\to\infty.
$$
\end{proposition}
\begin{proof}
The proof is based on the fact that for arbitrary $a>0$
$$
\int\limits_0^\infty t^{-1}J_2(at)\,dt=\frac{1}{2}.
$$
Indeed, according to (\ref{Resp}) and (\ref{u_1^delta}), we have
what
\begin{equation*}
r_N(t)= \frac{2}{t} J_2(2Nt),
\end{equation*}
and according to the previous identity
$$
\int\limits_0^{+\infty} r_N(t)\,dt =1,\quad N=1,2,\ldots.
$$
Substituting the test function $ \xi \in C ^ \infty (0, + \infty)
$  into the integral and replacing $ 2Nt = s $, we obtain
\begin{multline*}
\int\limits_0^{+\infty} r_N(t)\xi(t)\,dt=\int\limits_0^{+\infty}
\frac{2}{t} J_2(2Nt)\xi(t)\,dt=
\int\limits_0^{+\infty} \frac{2\cdot 2N}{s} J_2(s)\xi\left(\frac{s}{2N}\right)\,\frac{ds}{2N}=\\
=\int\limits_0^{+\infty} \frac{2}{s}
J_2(s)\xi\left(\frac{s}{2N}\right)\,ds \underset{N \to
\infty}{\longrightarrow} \int\limits_0^{+\infty} \frac{2}{s}
J_2(s)\xi(0)\,ds=(\delta(t),\xi(t)),
\end{multline*}
which completes the proof of the statement.
\end{proof}

The Krein-Stieltjes string with uniformly distributed identical
point masses is a discrete approximation of a unit density string.
For the initial-boundary value problem for a string with unit
density (\ref{eq_unit}) it is known that the response function has
the form $ r (t) = - \delta'(t) $. From Proposition \ref{p2} it
can be seen that the response function of the approximate system
does not converge to the response function of the original system.
This is due to the determination of the response function (the
reaction operator) in an approximate (discrete) system. Indeed,
for the reaction operator, we denoted an expression of the form $
R^Tf (t) = u_1^f (t) $ \, - \, the value of the solution in the
first channel (at the point $ x_1 = 1 / N $). In the original
system (\ref{eq_unit}), the reaction operator denotes the
expression $ R^Tf (t) = v_x^f (0, t) $ \, - \, the value of the
first derivative of the solution at the point $ x = 0 $ . If we
were to follow the ``correct'' reaction in a discrete system,
namely, an expression of the form
\begin{equation}\label{pr}
\tilde r_N=\frac{u_1^f(t)-u_0^f(t)}{1/N},
\end{equation}
 then we would prove the following
\begin{proposition}\label{p3}
Corrected response function $ \tilde r_N (t) $ \ref{pr} for
dynamic system (\ref{DynSyst}) with matrices bounded as (\ref{MN})
(for the Krein-Stieltjes string with uniformly distributed
identical point masses) converges in the sense of distributions to
the response function for the system (\ref{eq_unit})
$$
\tilde r_N(t)\to  r(t)=-\delta'(t), \quad\text{ if }N\to\infty.
$$
\end{proposition}

Using (\ref{con_v_f}) we conctruct a solution for (\ref{L_eq}).
Now we show that this solution converges in the sense of
distributions to $ v^\delta (x, t) $ \, - \, the solution of
(\ref{eq_unit}).
\begin{proposition}\label{p4}
The solution (\ref{con_v_f}) of the initial-boundary value problem
(\ref{L_eq}) converges in the sense of distributions
$$
u^\delta(x,t)\to v^\delta(x,t)=\delta(x-t)\quad\text{ if
}N\to\infty,
$$
where $v^\delta(x,t)$\,--\, the solution of the problem (\ref{eq_unit}).
\end{proposition}
\begin{proof}
We integrate the product of the function defined in
(\ref{con_v_f}) with arbitrary test function $ \xi \in C ^
\infty (0,1) $:
\begin{equation*}
(u^\delta(\cdot,t),\xi(\cdot))= \sum\limits_{j=1}\limits^{N}
(J_{2j-1}(2Nt)+J_{2j+1}(2Nt))\xi\left(\frac{j}{N}\right)+\underset{N
\to \infty}{o(1)}.
\end{equation*}
Taking $ \xi (x) = \sin (kx) $ as a test function, we see
that the right hand side of the last expression has the form
$$
2\sum\limits_{j=1}\limits^{N} J_{2j+1}(2Nt)
\sin\left(\frac{{(2j+1)k}}{2N}\right)+\underset{N \to
\infty}{o(1)}.
$$
Using the relation between the generating function and the
series \cite{SSF} (formula (9.1.43) on page 183):
$$
\sin(z\sin\theta)=2\sum\limits_{k=0}^{\infty}J_{2k+1}(z)\sin(2k+1)\theta,
$$
we come to the following expression
\begin{equation*}
\left(u^\delta(\cdot,t),\sin(k\cdot)\right)=
\sin\left(2Nt\sin\frac{k}{2N}\right)+\underset{N \to \infty}{o(1)}
=\sin(kt)+\underset{N \to \infty}{o(1)}.
\end{equation*}
Therefore,
\begin{equation*}
\left(u^\delta(\cdot,t),\sin(k\cdot)\right) \underset{N \to
\infty}{\longrightarrow}
\left(v^\delta(\cdot,t),\sin(k\cdot)\right) \quad \text{for all
$k\in \mathbb{N}$},
\end{equation*}
which means that there is convergence for an arbitrary test
function, which completes the proof.
\end{proof}

The statements \ref{p2} - \ref{p4} state that the solution of the
forward problem for a Krein-Stieltjes string with uniformly
distributed identical point masses, which is a discrete
approximation of a string with unit density, converges to the
solution of the original problem. A similar situation occurs in
the inverse problem \, --- \, the corrected response functions for
the discrete problem converge to the response functions of the
original one. In both forward and inverse problems, convergence
occurs in the space of distributions.

\section*{Acknowledgment}

The work of Victor Mikhaylov and Alexander Mikhaylov was supported
by  grants  Russian Foundation for Basic Research 18-01-00269. This work was supported by the program of the
Presidium of the Russian Academy of Sciences "The Modern Methods
of Mathematical Modeling in the Study of Nonlinear Dynamical
Systems" (targeted subsidy 08-04).

\end{document}